\documentclass[12pt,letterpaper]{article}

\title{Infinitely many 
conservation laws\\ for generalized nonlinear progressive wave equation}
\author{A. Sergyeyev\\[2mm]
{\small Mathematical Institute, Silesian University in Opava,}\\ {\small Na Rybn\'\i{}\v{c}ku 1, 74601 Opava, Czech Republic}\\ {\small E-mail \texttt{artur.sergyeyev@math.slu.cz}}}
\date{} 
\usepackage{amsmath, amssymb}
\newtheorem{theorem}{Theorem}
\begin{document}

\maketitle

\begin{abstract}
    We give a complete description of nontrivial local conservation laws of all orders for a natural generalization of the 
    nonlinear progressive wave equation and, in particular, show that there is an infinite number of such conservation laws. 

\end{abstract}
\section*{Introduction}
Consider a 
partial differential equation in $n+2$ independent variables $x,t,y_1,\dots,y_n$ and one dependent variable $u$,
\begin{equation}\label{gpe}
   u_{tx}+\left(f(u)\right)_{tt}+a u_{ttt}+b\Delta_y u=0,
\end{equation}
where $n$ is a nonnegative nonzero integer, 
$a$ and $b$ are nonzero constants, $f$ is an arbitrary 
nonlinear smooth function of $u$, and $\Delta_y=\sum_{i=1}^n\partial^2/\partial y_i^2$; the subscripts indicate partial derivatives in the usual manner.\looseness=-1

We shall hereinafter refer to this equation 
as to the {\em generalized nonlinear progressive wave equation} (GNPWE), as it is a natural generalization of the well-known nonlinear progressive wave equation (NP(W)E) 
which is recovered from (\ref{gpe}) if $f=c u^2$ for a constant $c$. 

Note that NPWE 
has numerous applications in the study of nonlinear wave processes in e.g.\  
meteorology, fluid dynamics, etc., see e.g.\ \cite{dkp, npe} and references therein. 
What is more, NPWE is, see e.g.\ \cite{dkp}, essentially the same as the KZK equation, which has important applications on its own right, for example in nonlinear acoustics.  

Our goal is to provide a complete description of 
all inequivalent nontrivial local conservation laws for (\ref{gpe}) with $a\neq 0$, $b\neq 0$ and smooth $f(u)$ such that $f_{uu}\neq 0$. 

Such a result is of interest as there are many important applications for conservation laws, cf.\ e.g.\ Ch.\ 4 in \cite{o} and references therein. For instance, discretizations of a PDE behave better when taking into account the (nontrivial) conservation laws admitted by the PDE in question, cf.\ e.g.\ \cite{bh}. 

On the other hand, finding 
all nontrivial local conservation laws for a given PDE is a difficult problem that was successfully addressed for only a rather small number of examples, see e.g.\ \cite{h,i, s} and references therein, and to the best of our knowledge GNPWE is not among those examples.

What is more, the overwhelming majority of the examples in question 
is strictly evolutionary, i.e., the PDEs under study 
have the form of the time derivative of the dependent variable being equal to a function of independent variables, dependent variable and its partial derivatives with respect to independent variables not involving time. GNPWE is not strictly evolutionary which makes the study of its conservation laws even more interesting.

\newpage
\section{Main Result}

\begin{theorem}\label{t}
For equation (\ref{gpe}) with any smooth function $f$ of $u$ such that $f_{uu}\neq 0$ and with any nonzero constants $a$ and $b$ all notrivial local conservation laws of all orders, which we write below as identities holding modulo (\ref{gpe}) and its differential consequences are, modulo the addition of trivial ones, of the following form:
\begin{equation}\label{ocl}
\bigl((u_x+f'(u)u_t+a u_{tt})\chi -(au_t+f(u))\chi_t\bigr)_{t}+(-u\chi_t)_x+\sum\limits_{j=1}^n \bigl(b\left(u_{y_j}\chi- u\chi_{y_j}\right)\bigr)_{y_j}=0,    
\end{equation}
where $\chi=\varphi_0+t\varphi_1$ and $\varphi_\alpha=\varphi_\alpha(x,y_1,\dots,y_n)$ are arbitrary smooth functions of their arguments satisfying the linear system 
\begin{equation}\label{sys}
\Delta_y\varphi_1=0,\quad \Delta_y\varphi_0=-(\varphi_1)_x/b.    
\end{equation}
\end{theorem}

Before proceeding to the proof note that the case of the original NP(W)E for which $f=cu^2$ with anonzero constant $c$ is {\em not} distinguished as far as the (local) conservation laws are concerned.

Also observe that for $n=1$ a general smooth solution of (\ref{sys}) takes a particularly simple form 
\[
\varphi_0=\eta_0+\eta_1 y_1-(\xi_0)_x y_1^2/(2b)-(\xi_1)_x y_1^3/(6 b),\quad 
\varphi_1=\xi_0+y_1\xi_1
\]
where $\eta_\alpha$ and $\xi_\alpha$, $\alpha=0,1$, are arbitrary smooth functions of $x$.

It is now clear that for any nonnegative nonzero integer $n$ equation (\ref{gpe}) 
has infinitely many nontrivial local conservation laws.

Note also that for $n=2,3,\dots$ it is possible to write down  a general smooth solution for (\ref{sys}) using integral representations for solutions of homogeneous and inhomogeneous Laplace equations but such representations are rather cumbersome and thus left beyond the scope of the present paper. 

\smallskip

\noindent{\em Proof.} 
Since equation (\ref{gpe}) is normal in the sense of \cite{o}, 
any nontrivial local conservation law for (\ref{gpe}) should have, see e.g.\ \cite{o}, a nonzero characteristic $\chi$ which is a differential function (recall that in our setting a differential function is a smooth function depending on $x,t,y_1,\dots,y_n,u$ and finitely many derivatives of $u$ and a local conservation law is a conservation law whose components are differential functions, cf.\ e.g.\ \cite{o}, to which we refer the reader for further details). 
Now observe that equation (\ref{gpe}) with arbitrary nonlinear smooth function $f$ of $u$ and nonzero constants $a$ and $b$ satisfies the conditions of Theorem 6 from \cite{i} with $O=-1$ 
and hence characteristics of 
local conservation laws of the 
equation under study may depend at most on the independent variables $x,t,y_1,\dots,y_n$.

Hence, any nontrivial local conservation law for (\ref{gpe}) is equivalent to a nontrivial local conservation law for which we have 
\begin{equation}
    \rho_{t} +\sigma_x +\sum\limits_{j=1}^n (\zeta_j)_{y_j}=\chi\cdot (u_{tx}+f_{tt}+a u_{ttt}+\Delta_y g),
\end{equation}
where the characteristic $\chi$ may depend at most on independent variables $x,t,y_1,\dots,y_n$ and $\rho,\sigma$, and $\zeta_j$, $j=1,\dots,n$, are local functions.

Then the characteristic $\chi=\chi(x,t,y_1,\dots,y_n)$ 
must satisfy the following equation, cf.\ e.g.\ Proposition 5.49 from \cite{o}: 
\begin{equation}\label{c} 
    \chi_{xt}+f_u\chi_{tt}-a\chi_{ttt}+b\Delta_y\chi=0
\end{equation}
It is a straightforward matter to verify that any smooth $\chi=\chi(x,t,y_1,\dots,y_n)$ satisfying (\ref{c}) is a characteristic for a local conservation law of the form (\ref{ocl}). 

To proceed, observe that as $f_{uu}$ is nonzero by assumption, $1$ and $f_u$ are linearly independent as functions of $u$.

Equating to zero the coefficients at linearly independent functions $f_u$ and 1 in (\ref{c}) yields
\begin{equation}\label{eqg}
    \chi_{tt}=0
\end{equation}
and 
\begin{equation}\label{eqa}
  \chi_{xt}-a\chi_{ttt} +b\Delta_y\chi=0. 
\end{equation}

From (\ref{eqg}) we get 
\[
\chi=\varphi_0+t\varphi_1
\]
where where $\varphi_\alpha=\varphi_\alpha(x,y_1,\dots,y_n)$, $\alpha=0,1$, 
are arbitrary smooth functions of their arguments 

Substituting this into (\ref{eqa}) and equating to zero the coefficients at the powers of $t$ yields
\[
(\varphi_1)_x+b\Delta_y\varphi_0=0
\]
and
\[
b\Delta_y\varphi_1=0,
\]
and, as $b\neq 0$ by assumption, the result follows. $\Box$

\subsection*{Acknowledgments}
The research of AS was supported in part through institutional support for the development of research organizations (RVO) for I\v{C} 47813059.

\end{document}